\documentclass[12pt,reqno,a4paper]{amsart}
\usepackage{amssymb,amscd, hyperref, epsf}

\begin{document}

\let\kappa=\varkappa
\let\epsilon=\varepsilon
\let\phi=\varphi
\let\p\partial

\def\Z{\mathbb Z}
\def\R{\mathbb R}
\def\C{\mathbb C}
\def\Q{\mathbb Q}
\def\P{\mathbb P}
\def\N{\mathbb N}
\def\L{\mathbb L}
\def\HH{\mathrm{H}}
\def\ss{X}

\def\conj{\overline}
\def\Beta{\mathrm{B}}
\def\const{\mathrm{const}}
\def\ov{\overline}
\def\wt{\widetilde}
\def\wh{\widehat}

\renewcommand{\Im}{\mathop{\mathrm{Im}}\nolimits}
\renewcommand{\Re}{\mathop{\mathrm{Re}}\nolimits}
\newcommand{\codim}{\mathop{\mathrm{codim}}\nolimits}
\newcommand{\id}{\mathop{\mathrm{id}}\nolimits}
\newcommand{\Aut}{\mathop{\mathrm{Aut}}\nolimits}
\newcommand{\lk}{\mathop{\mathrm{lk}}\nolimits}
\newcommand{\sign}{\mathop{\mathrm{sign}}\nolimits}
\newcommand{\pt}{\mathop{\mathrm{pt}}\nolimits}
\newcommand{\rk}{\mathop{\mathrm{rk}}\nolimits}
\newcommand{\SKY}{\mathop{\mathrm{SKY}}\nolimits}
\newcommand{\st}{\mathop{\mathrm{st}}\nolimits}
\def\Jet{{\mathcal J}}
\def\FC{{\mathrm{FCrit}}}
\def\sS{{\mathcal S}}
\def\lcan{\lambda_{\mathrm{can}}}
\def\ocan{\omega_{\mathrm{can}}}

\renewcommand{\mod}{\mathrel{\mathrm{mod}}}
\def\ds{\displaystyle}

\newtheorem{mainthm}{Theorem}
\newtheorem{thm}{Theorem}[section]
\newtheorem{lem}[thm]{Lemma}
\newtheorem{prop}[thm]{Proposition}
\newtheorem{cor}[thm]{Corollary}
\newtheorem{conjecture}[thm]{Conjecture}

\theoremstyle{definition}
\newtheorem{exm}[thm]{Example}
\newtheorem{rem}[thm]{Remark}
\newtheorem{df}[thm]{Definition}

\renewcommand{\thesubsection}{\arabic{subsection}}
\numberwithin{equation}{subsection}

\title[ Relations of Linking and Causality in Causally Simple Spacetimes]{Conjectures on the Relations of Linking and Causality in Causally Simple Spacetimes}
\author[Chernov]{Vladimir Chernov}
\thanks{This work was partially supported by a grant from the Simons Foundation
(\# 513272 to Vladimir Chernov). The author is thankful to Robert Low, Stefan Nemirovski and Miguel Sanchez for many useful discussions.}
\address{Department of Mathematics, 6188 Kemeny Hall,
Dartmouth College, Hanover, NH 03755, USA}
\email{Vladimir.Chernov@dartmouth.edu}

\begin{abstract} We formulate the generalization of the Legendrian Low conjecture of Natario and Tod (proved by Nemirovski and myself before) to the case of causally simple spacetimes. We prove a weakened version of the corresponding statement. 

In all known examples, a causally simple spacetime $(\ss, g)$ can be conformally embedded as an open subset into some globally hyperbolic $(\widetilde \ss, \widetilde g)$ and the space of light rays in $(\ss, g)$ is an open submanifold of the space of light rays in $(\widetilde \ss, \widetilde g)$. 
If this is always the case, this provides an approach to solving the conjectures relating causality and linking in causally simple spacetimes.

\end{abstract}

\maketitle

All manifolds, maps etc.~are assumed to be smooth unless the opposite is explicitly stated,
and the word {\it smooth\/} means $C^{\infty}$.

\section{Preliminaries and some definitions}\label{prelim}

A {\it spacetime\/} $(\ss, g)$ is a time oriented Lorentz manifold and points in it are called {\it events}. The {\it causal future\/} $J^+(x)$ of an event $x\in \ss$ is the set of all points in $\ss$ that one can reach from $x$ going into the future with less or equal than light speed. The {\it chronological future\/} $I^+(x)$ of $x\in \ss$ is the set of all points that one can reach from $x$ going into the future with less than light speed. Causal past $J^-(x)$ and chronological past $I^-(x)$  are defined similarly.

Two Lorentz manifolds $(X_1, g_1)$ and $(X_2, g_2)$ are said to be {\em conformal\/} if there exists
a diffeomorphism $\Phi:X_1\to X_2$ and a positive smooth function $\Omega:X_1\to(0, +\infty)$ 
such that $\Omega g_1=\Phi^*(g_2).$ 

Clearly for all $x\in X_1,$ we have that 
$\Phi(J^{\pm}(x))=J^{\pm}(\Phi(x))$ and $\Phi(I^{\pm}(x))=I^{\pm}(\Phi(x))$. So in particular
if $(X_1, g_1)$ is causally simple, or strongly causal, or globally hyperbolic, then 
$(X_2, g_2)$ is respectively causally simple, or strongly causal, or globally hyperbolic

Moreover if $\gamma(t)$ is a null pregeodesic for $(X_1, g_1)$, 
then $f(\gamma(t))$ is a null pregeodesic for $(X_2, g_2)$, see~\cite[Lemma 9.17]{BeemEhrlichEasley}.
Such a statement is generally false for timelike and spacelike pregeodesics.

A spacetime is {\it globally hyperbolic\/} if it does not have closed trajectories representing movement with less or equal than light speed and for any two events $x,y$ the intersection 
$J^+(x)\cap J^-(y)$ is compact. (See~\cite{BernalSanchezGloballyHyperbolic} for the proof that this definition of a globally hyperbolic spacetime is equivalent to the standard one.) 

A spacetime $(X^{m+1}, g)$ is called {\it causally simple\/} if  if it does not have closed trajectories representing movement with less or equal than light speed {\bf and\/} for every $x\in X$ the sets $J^+(x)$ and 
$J^-(x)$ are closed. This class of spacetimes is important because it is the best behaved class in the ladder of causality~\cite{MinguzziSanchez} that allows naked singularities.

\section{Results}

We start by proving the following Theorem.

\begin{thm}\label{WeakLow}
Let $(\ss^{m+1}, g), m>1$ be a causally simple spacetime and let $S^0=\{\pt_1, \pt_2\}$ be the two point 
space. 
Let $\phi:S^0\times [0,1]\to X$ be a continuous map, such that the image of $\phi|_{S^0\times 0}:S^0\to X$
is two causally related points and the image of $\phi|_{S^0\times 1}:S^0\to X$
is two causally unrelated points. Then there is $\tau\in [0,1]$ such that 
$\phi(\pt_1, \tau)$ and $\phi(\pt_2, \tau)$ are located on the same null geodesic.
\end{thm}

Implication ${\bf 2}\implies {\bf 3}$ in~\cite[Theorem 15]{ChernovRudyakCMP} says that the statement of
Theorem~\ref{WeakLow} is true for all globally hyperbolic spacetimes. As it is explained 
in~\cite[Remark 9]{ChernovRudyakCMP}, the proof can be generalized to all 
causally simple spacetimes for which the Lorentz distance function is a finite continuous function. 
Recent results of Mueller and Sanchez~\cite{MuellerSanchez} allow us to prove this result for 
all causally simple spacetimes.

\begin{proof}
Let us recall the
notion of Lorentzian distance, see~\cite{BeemEhrlichEasley}. For
two points $p,q$ in $(X,g)$ with $q\in J^+(p)$ put $\Omega_{p,q}$ to be the space of all
piecewise smooth future directed non-spacelike curves
$\alpha:[0,1]\to X$ with $\alpha(0)=p, \alpha(1)=q.$ For $\alpha\in
\Omega_{p,q}$ choose a partition $0=t_0<t_1<t_2\cdots
<t_{k-1}<t_k=1$ such that $\alpha |_{(t_i, t_{i+1})}$ is smooth for
all $i\in\{0,1,\cdots(k-1)\},$ and define the {\em Lorentzian arc
length\/} $L(\alpha)$ of $\alpha$ by
$$
L(\alpha)=L_g(\alpha)=\sum_{i=0}^{k-1}\int_{t_i}^{t_{i+1}}\sqrt{-g(\dot
{\alpha}(s), \dot {\alpha}(\tau))}d s.
$$

For $p,q\in (X,g)$ define the {\em Lorentzian distance function \/}
$d=d_g:X\times X\to \R\sqcup\infty$ as follows: set $d(p,q)=0$ for
$q\not \in J^+(p);$ and set $d(p,q)=\sup\{L_g(\alpha)|\alpha\in
\Omega_{p,q}\},$ for $q\in J^+(p).$

By~\cite[Chapter 14, Corollary
1]{ONeill}, if $b\in I^+(a)$ and $c\in J^+(b)$ or if $b\in J^+(a)$ and
$c\in I^+(b),$ then $c\in I^+(a).$ Combining this with the definition of
the Lorentzian distance $d$ we get that $d(p,q)>0$ if and only if
$q\in I^+(p)$.

In general, the Lorentzian distance function is not continuous, and
$d(p,q)$ is not finite. 

Every causally simple spacetime is stably causal, see 
for example~\cite[Proposition 3.61 and Remark 3.66]{MinguzziSanchez}. 
Mueller and Sanchez~\cite[Corollary 1.4]{MuellerSanchez} proved 
that since $(X,g)$ is stably causal, it admits a conformal embedding $\Psi$ into the Minkowski spacetime 
$(\R^{N+1}, g_{\st}=\sum_{i=1}^Ndx_i^2-dt^2)$ of large dimension. 
Since $(\R^{N+1}, g_{\st})$ is globally hyperbolic, the Lorentz distance function on it
is a continuous function with values in $\R,$ rather than in $\R\sqcup \infty$, 
see~\cite[Corollary 4.7]{BeemEhrlichEasley}. 
Thus the Lorentz distance function on $X$ defined by 
$\Psi^*(g_{\st})$  takes values in $\R$, but it is not necessarily continuous. 

Since 
$(X, \Psi^*(g_{\st}))$ is conformal to the causally simple $(X,g)$, it also is causally simple.
Clearly the statement of Theorem~\ref{WeakLow} is true for $(X, g)$ if and only 
if it is true for $(X, \Psi^*(g_{\st}))$. Thus it suffices to prove the Theorem for a causally simple $(X,g)$
such that the induced Lorentz distance is a (not necessarily continuous) function with values in $\R$, rather 
than in $\R\sqcup \infty$.
 
Put $y_i=\phi(\pt_i, 0)$ and $x_i=\phi(\pt_i,1), i=1,2.$
{\em We argue by contradiction and assume that there is no $\tau\in [0,1]$ such that 
$\phi(\pt_1, \tau)$ and $\phi(\pt_2, \tau)$ are on the same null geodesic.\/}
Since $y_1, y_2\in \ss$ are
causally related, either $y_1\in J^+(y_2)$ or $y_2\in J^+(y_1).$
Without loss of generality assume that $y_2\in J^+(y_1).$ If $y_2\in
J^+(y_1)\setminus I^+(y_1),$ then $y_1=\phi(\pt_1, 0)$ and $y_2=\phi(\pt_2, 0)$ lie on a common
null geodesic, see~\cite[Corollary 4.14]{BeemEhrlichEasley}. This contradicts to our assumptions.
Hence $y_2\in I^+(y_1)$ and $d(y_1, y_2)>0.$

Define the function $\ov d:[0,1]\to \R$ by $\ov
d(t)=d(\phi(\pt_1, t), \phi(\pt_2, t)).$ This function generally is not continuous. 
We have $\ov d(0)=d(p_1(0), p_2(0))=d(y_1,
y_2)>0.$ Furthermore, $\ov d(1)=d(p_1(1), p_2(1))=d(x_1, x_2)=0,$
since $x_1, x_2$ are causally unrelated. Put $\tau=\inf \{t\in
[0,1]| \ov d(t)=0\}$, so that
\begin{equation}\label{positivedeq1}
\ov d(t)>0 \text{ for all } t<\tau.
\end{equation}

By~\cite[Lemma 4.4]{BeemEhrlichEasley} we have that if $d(p,q)<\infty$, $p_n\to p$ and $q_n\to q$, then 
$d(p,q)<\lim \inf d(p_n, q_n).$ Taking $p=\phi(\pt_1, \tau), q=\phi(\pt_2, \tau)$ we 
get that
\begin{equation}\label{positivedeq2}
\ov d(\tau)=0.
\end{equation}

Below we show that $\phi(\pt_1, \tau)$ and $\phi(\pt_2, \tau)$ belong to the same null geodesic, thus getting 
a contradiction. Starting from here we pretty much repeat the corresponding part of the proof of 
implication ${\bf 2}\implies {\bf 3}$ in~\cite[Theorem 15]{ChernovRudyakCMP}.

By~\cite[Proposition 3.68]{MinguzziSanchez} since $(X,g)$ is causally simple,
the sets $J^{\pm}(K)=\cup_{k\in K}J^{\pm}(k)$ are
closed for every compact $K\subset \ss$.

Define the continuous path $p_1, p_2:[0,1]\to X$ by $p_1(t)=\phi(\pt_1, t), p_2(t)=\phi(\pt_2, t)$. 
By \eqref{positivedeq1}, $d(p_1(t), p_2(t))=\ov d(t)>0$ for all
$t<\tau.$ Thus $p_2(t)\in I^+(p_1(t))\subset J^+(p_1, t))$ for every
$t<\tau,$ and so $\Im (p_2|_{[t, \tau)})\subset J^+(\Im (p_1|_{[t,
\tau]}))$ for every $t<\tau.$ Since $\Im (p_1|_{[t, \tau]})$ is
compact and $(\ss, g)$ is causally simple, we conclude that $J^+(\Im
(p_1|_{[t, \tau]}))$ is closed, and thus $p_2(\tau)\in J^+(\Im
(p_1|_{[t, \tau]}))$ for every $t<\tau.$

We choose an increasing sequence $\{t_i\in [0,1]\}_{i\in \N}$ that
converges to $\tau.$ Then for each $i\in \N$ there exists $\wt
t_i\in [t_i, \tau]$ such that $p_2(\tau)\in J^+(p_1(\wt t_i)).$
Thus $p_1(\wt t_i)\in J^-(p_2(\tau))$ for all $i.$ Since $(\ss, g)$
is causally simple, $J^-(p_2(\tau))$ is closed and it contains the
point $\ds p_1(\tau)=\lim_{i\to \infty} p_1(\wt t_i).$ Since
$p_1(\tau)\in J^-(p_2(\tau)),$ we have $p_2(\tau)\in
J^+(p_1(\tau)).$

On the other hand, $p_2(\tau)\not \in I^+(p_1(\tau))$ since
$d(p_1(\tau), p_2(\tau))=\ov d(\tau)=0$ by \eqref{positivedeq2}. So,
$p_2(\tau)\in J^+(p_1(\tau))\setminus I^+(p_1(\tau))$, and therefore
the points $p_1(\tau), p_2(\tau)$ belong to a common null geodesic,
see~\cite[Corollary 4.14]{BeemEhrlichEasley}. 
\end{proof}

\begin{rem}[importance of being simple]
The statement of Theorem~\ref{WeakLow} is generally false for spacetimes that are not causally simple. Indeed 
put $(X^{1+1}, g)$ to be the spacetime obtained by deleting the point $\{(1,1)\}$ from
the Minkowski spacetime $(\R^2, dx^2-dt^2).$ Define the continuous $\phi:S^0\times [0,1]\to X$  by 
$\phi(\pt_1, t)=(0,0)$ and $\phi(\pt_2, t)=(4t, 2)$, $t\in [0,1]$. 
The points $\{(0,0), (0,2)\}=\{\phi(\pt_1, 0), \phi(pt_2, 0)\}$ are causally related, the points
$\{(0,0), (4,2)\}=\{\phi(\pt_1, 1), \phi(pt_2, 1)\}$ are causally unrelated. 
Note that for $\tau=\frac{1}{2}$ the points
$\phi(\pt_1, \tau)$ and  $\phi(\pt_2, \tau)$ are located on the same null geodesic in the Minkowski 
spacetime, but there is no $\tau$  such that $\phi(\pt_1, \tau)$ and $\phi(\pt_2, \tau)$ are on the same 
null geodesic in $(X^{1+1}, g)$.
\end{rem}

A seminal observation of Penrose and Low~\cite{PR2, Low0, LowNullGeodesics} is that the space $\mathcal N_X$ of all the unparameterized future directed light rays in (strongly causal) $\ss$
has a canonical structure of a (possibly non-Hausdorff) contact manifold (see also \cite{NatarioTod, KhesinTabachnikov, BIL}). For a globally hyperbolic spacetime the ratio of the contact forms on the space of light rays coming from two Cauchy surfaces equals to the redshift between them~\cite{ChernovNemirovskiRedshift}.

\begin{rem} 
The null geodesics through $x\in X$ give an embedded $(m-1)$-sphere $\mathfrak S_x\subset \mathcal N_X$
called the {\it sky of $x$.} 
The sky $\mathfrak S_x$ is a Legendrian submanifold of $\mathcal N_X.$ Each causally simple 
spacetime is strongly causal, see for example E.~Minguzzi and 
M.~Sanchez~\cite[Proposition 3.57, Proposition 3.61 and Remark 3.66]{MinguzziSanchez}.

It is reasonable to ask whether given 
any two pairs of causally unrelated points $\{x_1, x_2\}$ and $\{y_1, y_2\}$ in a causally simple spacetime 
$(X^{m+1}, g),$ $m>1$ there always exists a continuous map $\phi:S^0\times [0,1]\to X$ such that 
$\phi(\pt_i, 0)=x_i$, $\phi(\pt_i, 1)=y_i, i=1,2$ and for all $t\in [0,1]$ the two points
$\phi(\pt_1, t)$ and $\phi(\pt_2, t)$ do not belong to the same null geodesic. 

If this statement is true then for a causally simple spacetime $(X,g)$ 
the Legendrian links $(\mathfrak S_{x_1}, \mathfrak S_{x_2})$ and $(\mathfrak S_{y_1}, \mathfrak S_{y_2})
\subset \mathcal N_X$ formed by skies of 
two pairs of causally unrelated points are always Legendrian isotopic, with the isotopy given by 
$\bigl(\mathfrak S_{\phi(\pt_1, t)}, \mathfrak S_{\phi(\pt_2, t)}\bigr), t\in [0,1].$ This would show that any two two-component links consisting of skies of two causally unrelated events in a causally simple spacetime are Legendrian isotopic. The class of such a link would be a natural candidate for the {\it unlink.\/}
\end{rem}

\section{Legendrian Low conjecture for nonrefocussing causally simple spacetimes and some other conjectures.\/}\label{s:refocussing}
Robert Low~\cite{LowNullGeodesics} introduced the following definition of nonrefocussing spacetimes.

\begin{df}\label{d:nonrefocussing}
A strongly causal spacetime $(\ss,g)$ is (weakly) {\em refocussing at $x\in X$\/} if
there exists a neighborhood $V\ni x$ satisfying the following property:
For every open $U$ with $x\in U \subset V$ there exists $y\not \in
U$ such that all the null-geodesics through $y$ enter $U.$ A
spacetime $(X,g)$ is called {\em refocussing\/} if it is
refocussing at some $x,$ and it is called {\em nonrefocussing\/} if
it is not refocussing at every $x\in X$.

A strongly causal spacetime $(\ss,g)$ is {\em strongly refocussing at $x\in X$\/} if there exists 
$y\neq x$ such that all null geodesics through $y$ pass through $x.$ A strongly causal spacetime is 
{\em strongly refocussing\/} if it is strongly refocussing at some $x$.
\end{df}

\begin{rem}[Refocussing versus strong refocussing]
Clearly every strongly refocussing spacetime is refocussing. At the moment we do not know
any example of spacetimes that are refocussing but not strongly so. 
However Kinlaw~\cite{PaulKinlaw} and~\cite[Example 5.6]{CKS} has constructed examples of 
spacetimes that are refocussing at some $x$ but not strongly refocussing at $x$. 

The simplest example is the globally hyperbolic Einstein cylinder 
$(S^n\times (-\pi, \pi), g_1\oplus -dt^2)$, where $g_1$ is the standard Riemannian metric on the unit sphere 
$S^n\subset \R^{n+1}.$ This spacetime is refocussing but not strongly refocussing at each $(x,0).$ However 
it is strictly refocussing at each $(x,t), t\neq 0.$
Kinlaw's examples can be made null and timelike geodesically complete, however they all are strongly 
refocussing at some other $\hat x$.

Also in his Ph.D. Thesis Levi~\cite{Levi} proved that every refocussing two-dimensional spacetime is strongly refocussing.
\end{rem}

Let $\SKY$ denote the set of all skies in $(\ss, g)$ with the following
topology. For $\mathfrak S_x\in \SKY$ the topology base at $\mathfrak S_x$ is given 
by $\{\mathfrak S_y| \mathfrak S_y\subset W\}$ for open $W\subset \mathcal N$ such 
that $\mathfrak S_x\subset W.$

Consider the map
\begin{equation}\label{eq:mu}
\mu: \ss\to \SKY, \quad \mu(x)=\mathfrak S_x.
\end{equation}

Low~\cite{LowNullGeodesics} observed that if a strongly causal $(X^{m+1}, g)$ is nonrefocussing, then 
$\mu$ is a homeomorphism (see also Bautista, Ibort, Lafuente~\cite{BIL2, BIL3} for the recent development of Low's work). The proof of this statement of Low
is contained in the PhD thesis of Kinlaw~\cite{PaulKinlaw}. Thus if $\rho:[0,1]\to \SKY$ is a
continuous path, then there is a continuous path $p:[0,1]\to X$ such that 
$\rho(t)=\mathfrak S_{p(t)}.$ 

The fact that points $q_1, q_2\in X$ are on the same null geodesic means exactly that 
$\mathfrak S_{q_1}\cap \mathfrak S_{q_2}\neq \emptyset,$ i.e. that the Legendrian 
link $(\mathfrak S_{q_1}, \mathfrak S_{q_2})\subset \mathcal N_X$ is singular.

Combining this with Theorem~\ref{WeakLow} we immediately get the following Corollary.

\begin{cor}\label{WeakLowCorollary}
Let $(X^{m+1}, g), m>1$ be a causally simple nonrefocussing spacetime. Let $x_1, x_2\in X$ 
be two causally unrelated points, and let $y_1, y_2\in X$ be two causally related points 
that do not lie on the same null geodesic, so that the Legendrian link 
$(\mathfrak S_{y_1},\mathfrak S_{y_2})\subset \mathcal N_X$ is nonsingular. Then there is 
no Legendrian isotopy changing the link $(\mathfrak S_{y_1},\mathfrak S_{y_2})$ into 
$(\mathfrak S_{x_1},\mathfrak S_{x_2})$ that satisfies the additional restriction that at each moment of the 
isotopy the corresponding Legendrian link consists of skies of points in $(X,g).$
\end{cor}

This motivates the {\it ``Legendrian Low conjecture'' for causally simple spacetimes\/} that we formulate below. 

\begin{conjecture}\label{LowConjectureCausallySimple}
Let $(X^{m+1},g),$ $m>1$ be a causally simple spacetime 
that does not admit a metric $g'$ such that $(X, g')$ is refocussing and causally simple. 
Then the
Legendrian links formed by a pair of skies of two causally unrelated events in $(X,g)$ and by a pair of skies 
of two causally related events in $(X,g)$ 
(that do not belong to the same null geodesic) are never Legendrian isotopic.
\end{conjecture}

Clearly Corollary~\ref{WeakLowCorollary} is saying that a certain weakened version of the above
conjecture is indeed true. This conjecture is motivated by the following discussion.

For globally hyperbolic spacetimes the Legendrian links formed by skies of any two pairs of 
causally unrelated points are Legendrian isotopic, see for example~\cite[Theorem 8]{ChernovRudyakCMP} 
and~\cite[Lemma 4.3]{ChernovNemirovskiGAFA}. The Legendrian isotopy class of such a link is called 
the class of the {\it ``unlink''.\/}

The {\em Low conjecture}~\cite{Low0} asked whether it is true that the link formed by skies of two causally 
related events in $(X^{2+1}, g)$ is never topologically isotopic to the unlink, provided that 
$(X^{2+1}, g)$ is globally hyperbolic and its Cauchy surface $M^2$ is homeomorphic to a disk with holes.

The {\em Legendrian Low conjecture\/} formulated by Natario and Tod~\cite{NatarioTod} 
asked whether it is true that the Legendrian link formed by skies of two causally related 
events in $(X^{3+1}, g)$ is never Legendrian isotopic to the unlink, provided that 
$(X^{3+1}, g)$ is globally hyperbolic and its Cauchy surface $M^3$ is homeomorphic to $\R^3.$

Jointly with Nemirovski we proved these conjectures in~\cite[Theorem A and Theorem B]{ChernovNemirovskiGAFA}. 
In fact we proved that 
the statement of the Low conjecture remains true 
whenever the universal cover of the two-dimensional Cauchy surface $M^2$  of $(X^{2+1}, g)$ 
is not a closed manifold. Our result~\cite[Theorem 10.4]{ChernovNemirovskiGT} 
says that the statement of the Legendrian Low conjecture is true for all globally 
hyperbolic spacetimes  $(X^{m+1}, g), m>1$ provided that the universal cover of its $m$-dimensional Cauchy 
surface is not a closed  manifold. See also~\cite{ChernovNemirovskiOrder, ChernovCompactUniversalCovering} for the generalizations of these results.

If $(X^{m+1}, g)$ is refocussing, then the universal covering of its Cauchy surface is a
closed manifold, see~\cite{LowNullGeodesics, LowRefocussing} 
and~\cite[Proposition 6 and Theorem 14]{ChernovRudyakCMP}. 
Similarly to~\cite[Example 10.5]{ChernovNemirovskiGT} one can show that  the Legendrian Low 
conjecture fails for every strongly refocussing globally hyperbolic $(X^{m+1}, g)$, $m>1.$ 

When $m=2,3$ and the universal cover of a Cauchy surface $M$ is a closed manifold, it has to be homeomorphic to the sphere $S^m$. Moreover one can equip the Cauchy surface with the quotient $\overline g$ of the round sphere metric. This is trivial for $m=2$, while for $m=3$ this is the Poincare conjecture and the Thurston's Elliptization Conjecture proved by Perelman~\cite{Perelman1, Perelman2, Perelman3}. Clearly the globally hyperbolic spacetime $(M\times \R, \overline g\oplus -dt^2)$ is strongly refocussing. 

Thus we get that causality is equivalent to Legendrian linking of skies  for globally hyperbolic spacetimes $X^{m+1}, m=2,3$ for which one can not put a refocussing metric on $X$. We think that the same fact is true in higher dimensions, namely:

\begin{conjecture}\label{LowConjectureNonrefocussing} 
Given a globally hyperbolic $(X^{m+1}, g), m \geq 2,$ causal relation of two points is equivalent to Legendrian linking of their skies, provided that one can not put a globally hyperbolic and refocussing metric on $X$.
\end{conjecture}

So essentially our Conjecture~\ref{LowConjectureCausallySimple} generalizes the above Conjecture~\ref{LowConjectureNonrefocussing} to the case of causally simple spacetimes.

Below we formulate two more conjectures that, if they are true, should give an approach to solving Conjecture~\ref{LowConjectureCausallySimple}.

In all known examples of a causally simple spacetime $(\ss, g)$ it can be conformally embedded into a globally hyperbolic spacetime $(\widetilde \ss, \widetilde g)$ of the same dimension. So we conjecture that this is always the case.

\begin{conjecture}\label{conformalembedding}
Let $(\ss, g)$ be a causally simple spacetime, then it can be conformally embedded as an open subset into some globally hyperbolic $(\widetilde \ss, \widetilde g).$
\end{conjecture}

Moreover in all the known examples where a causally simple $(\ss, g)$ is conformally embedded as an open set into a globally hyperbolic $(\widetilde \ss, \widetilde g)$, the space of (future directed, unparameterized) light rays in $(\ss, g)$ appears to be an open submanifold of the space of light rays of $(\widetilde \ss, \widetilde g)$. So we conjecture that this is always the case.

\begin{conjecture}\label{LightRaySubspace}
Let $(\ss, g)$ be a causally simple spacetime that can be conformally embedded as an open subset into some globally hyperbolic $(\widetilde \ss, \widetilde g)$, then the space of light rays in $(\ss, g)$ is an open contact submanifold of the space of light rays in $(\widetilde \ss, \widetilde g).$  
\end{conjecture}

Since a light ray in $(\widetilde \ss, \widetilde g)$ can be split into many light rays in $(\ss, g)$, we do not see why Conjecture~\ref{conformalembedding} would immediately imply Conjecture~\ref{LightRaySubspace} or even that the space of light rays in $(\widetilde \ss, \widetilde g)$ is Hausdorff.

We finish by showing that Conjecture~\ref{LowConjectureCausallySimple} is true for causally simple $(\ss,g)$  satisfying Conjecture~\ref{LightRaySubspace} and the following two assumptions:
\begin{description}
\item[1] in $(\widetilde \ss, \widetilde g)$ we know that the Legendrian link consisting of skies of any two causally related points is not Legendrian isotopic to the unlink.
(So for example, we could take any $(\widetilde \ss, \widetilde g)$ for which the universal cover of its Cauchy surface is not compact~\cite{ChernovNemirovskiGT}; or for which the universal cover of its Cauchy surface is compact, but the integral cohomology ring of the universal cover is not the one of a CROSS  (compact rank  one symmetric space) \cite{ChernovCompactUniversalCovering}.)
\item[2] every pair of causally unrelated events $x,y\in \ss$ remains to be causally unrelated when considered to be a pair of events in $(\widetilde \ss, \widetilde g).$
\end{description}

We argue by contradiction and suppose that $x,y\in \ss$ are causally related, $\hat x, \hat y\in \ss$ are causally unrelated, but the Legendrian links 
$(\mathfrak S_x, \mathfrak S_y)$ and $(\mathfrak S_{\hat x}, \mathfrak S_{\hat y})$ are Legendrian isotopic in the contact manifold $\mathcal N_{\ss}$ of all light rays in $\ss.$ Then $\hat x, \hat y$ are causally related considered as points of $(\widetilde \ss, \widetilde g),$ while by our assumption $x,y,$ are causally unrelated in $(\widetilde \ss, \widetilde g)$; and clearly the Legendrian links 
$(\mathfrak S_x, \mathfrak S_y)$ and $(\mathfrak S_{\hat x}, \mathfrak S_{\hat y})$ are Legendrian isotopic in the contact manifold $\mathcal N_{\widetilde \ss}$ of all light rays in $\widetilde \ss.$ However for $\widetilde \ss$ satisfying condition $1$ above this is not possible, so we got a contradiction.

\end{document}